\documentclass[12pt]{amsart}
\usepackage{newlfont,amsmath,enumerate,amssymb}

\newcommand{\overgamma}{\overline{\gamma}}
\newcommand{\overzeta}{\overline{\zeta}}
\newcommand{\phione}{F_{-1}}
\newcommand{\phitwo}{F_2}
\newcommand{\phij}{F_j}

\newcommand{\calH}{{\mathcal{H}}}

\newcommand{\calR}{{\mathcal{R}}}

\newcommand{\PGL}{{\mathrm{PGL}}}

\newcommand{\rmss}{{\mathrm{ss}}}

\newcommand{\bbA}{{\mathbb{A}}}
\newcommand{\bbC}{{\mathbb{C}}}

\newcommand{\bbQ}{{\mathbb{Q}}}
\newcommand{\bbR}{{\mathbb{R}}}
\newcommand{\bbZ}{{\mathbb{Z}}}

\newcommand{\rSh}{{\mathrm{Sh}}}

\newcommand{\pr}{{\mathrm{pr}}}

\newcommand{\End}{{\mathrm{End}}}
\newcommand{\Ind}{{\mathrm{Ind}}}

\newcommand{\SL}{{\mathrm{SL}}}
\newcommand{\GL}{{\mathrm{GL}}}

\newcommand{\bfq}{\mathbf{q}}

\newcommand{\bfs}{\mathbf{s}}

\newcommand{\bfT}{\mathbf{T}}
\newcommand{\bfW}{\mathbf{W}}

\newtheorem{lemma}{Lemma}

\newtheorem{prop}[lemma]{Proposition}
\newtheorem{thm}[lemma]{Theorem}
\newtheorem{cor}[lemma]{Corollary}

\begin{document}


\title[Spherical vectors]{A representation theoretic approach to
   Kohnen's plus space of modular forms of half integral weight}

\author{Hung Yean Loke and Gordan Savin}

\address{Hung Yean Loke, Department of Mathematics, National
  University of Singapore, 2 Science Drive 2, Singapore 117543}
\email{matlhy@nus.edu.sg} \address{Gordan Savin, Department of
  Mathematics, University of Utah, Salt Lake City, UT 84112}
\email{savin@math.utah.edu}

\begin{abstract}
In this paper, we define a notion of pseudo-spherical type for the two 
fold central extension of $\SL_2(\bbQ_2)$. We relate this definition to some
results in classical modular forms of half integral weights.
\end{abstract}

\maketitle

\section{Introduction} 

Let $\bbQ$ be the field of rational numbers. For every place $v$ of
$\bbQ$ let $\bbQ_v$ denote the corresponding local field.  Then
$\bbQ_v=\bbR$ or $\bbQ_p$ for a prime $p$.  The group $\SL_2(\bbQ_v)$
has a non-trivial two-fold central extension
\begin{equation} \label{eq1} 
1\rightarrow \mu_2 \rightarrow G(\bbQ_v) \rightarrow \SL_2(\bbQ_v)
\rightarrow 1
\end{equation}
where $\mu_2 = \{ \pm 1 \}$.  Recall that an irreducible
representation of $G(\bbQ_{v})$ is called genuine if the central
subgroup $\mu_{2}$ acts faithfully on it. Gelbart's book \cite{G2}
contains a basic theory of genuine representations of $G(\bbR)$ and
$G(\bbQ_{p})$ for $p\neq 2$. Our intent is to develop a theory in the
case of $G(\bbQ_{2})$. The main difference between $G(\bbQ_{2})$ and
$G(\bbQ_{p})$ for $p\neq 2$ lies in the fact that the central
extension splits over $\SL_{2}(\bbZ_{p})$ for $p\neq 2$. In
particular, we have a subgroup $K_{p}\subseteq G(\bbQ_{p})$ isomorphic
to $\SL_{2}(\bbZ_{p})$ under the natural projection from $G(\bbQ_{p})$
to $\SL_{2}(\bbQ_{p})$ for every $p\neq 2$. A genuine representation
$\pi$ of $G(\bbQ_{p})$ is called {\it unramified} if it contains a
non-zero $K_{p}$-fixed vector.

Assume now that $p=2$. Let $K$ denote the full inverse image of
$\SL_{2}(\bbZ_{2})$ in $G(\bbQ_{2})$. In this case the central
extension splits over a smaller subgroup.  More precisely, we have a
subgroup $K_{1}(4)\subseteq K$ isomorphic to the subgroup of
$\SL_{2}(\bbZ_{2})$ given by the following congruence:
\[
\begin{pmatrix}
a & b \\
c & d
\end{pmatrix}
\equiv 
\begin{pmatrix}
1 & \ast \\
0 & 1
\end{pmatrix}
\pmod{4}
\]
In this paper we completely describe genuine irreducible
representations of $G(\bbQ_{2})$ containing non-zero $K_{1}(4)$-fixed
vectors.  More precisely, in Section \ref{S3}, we describe a Hecke algebra 
$H(\gamma)$ which captures the structure of all representations generated by 
$K_{1}(4)$-fixed vectors and with a fixed central character $\gamma$. 
In Section \ref{S4} we show that $H(\gamma)$ is isomorphic to the
Iwahori-Matsumoto Hecke algebra for $\PGL_{2}(\bbQ_{2})$. In this way we get 
a correspondence between (some) representations of $G(\bbQ_{2})$ and
representations of $\PGL_{2}(\bbQ_{2})$.
 We call this correspondence a local Shimura correspondence.  

In Section \ref{S5} we show that the compact group $K$ has exactly two
irreducible genuine representations, with the fixed central character
$\gamma$, containing non-zero $K_{1}(4)$-fixed vectors. These
representations are denoted by $V(2)$ and $V(-1)$ and have dimensions
2 and 4, respectively.  We show that a representation $\pi$ of $G(\bbQ_{2})$ has
$V(2)$ as a $K$-type if and only if it corresponds to an unramified
representation of $\PGL_{2}(\bbQ_{2})$, by the local Shimura 
correspondence.  Thus, it is natural to define unramified 
representations of $G(\bbQ_{2})$ to be those that contain $V(2)$ as a $K$-type, and 
we call $V(2)$ a pseudo-spherical type.

We should point out that the center of $G(\bbQ_{2})$ is a cyclic group
of order $4$. Thus, we have two different genuine central characters $\gamma$
and two classes of unramified representations. This is analogous to
the case of the real group $G(\bbR)$, where the weights $-1/2$ and
$1/2$ are called pseudo-spherical types.

We then apply our local results in a global setting in Section
\ref{S8}. Let $\bbA$ be the ring of adeles and let $G(\bbA)$ be the
two-fold cover of $\SL_{2}(\bbA)$.  Let $r>1$ be an odd integer. Let
$\pi=\otimes\pi_{v}$ be a genuine cuspidal automorphic representation
such that
\begin{itemize} 
\item $\pi_{\infty}$ is a holomorphic discrete series representation
  with the lowest weight $r/2$.
\item $\pi_{p}$ is unramified for all $p\neq 2$. 
\item $\pi_{2}$ contains a non-zero $K_{1}(4)$-fixed vector. 
\end{itemize} 
Every such $\pi$ corresponds to a Hecke eigenspace in
$S_{r/2}(\Gamma_{0}(4))$, the space of cuspidal modular forms of
weight $r/2$.  Roughly speaking, a function $f=\otimes f_{v}$ in $\pi$
such that $f_{\infty}$ is a lowest weight vector in $\pi_{\infty}$,
$f_{p}$ is $K_{p}$-fixed and $f_{2}$ is $K_{1}(4)$-fixed, gives
naturally a modular form in $S_{r/2}(\Gamma_{0}(4))$. Since the space
of $K_{1}(4)$-fixed vectors in $\pi_{2}$ is 2-dimensional, unless
$\pi_{2}$ is a Steinberg representation, the cuspidal automorphic
representation $\pi$ gives rise to a two-dimensional Hecke eigenspace
in $S_{r/2}(\Gamma_{0}(4))$. We can pick a line in this subspace by
taking $f_{2}$ to be in the $K$-type isomorphic to $V(2)$.  In this
way we get a representation-theoretic description of Kohnen's plus
space $S^{+}_{r/2}(\Gamma_{0}(4))$ \cite{Ko}. Finally, we show that the
global Shimura correspondence is compatible with our local Shimura
correspondence at the place $p=2$.

\section{Double cover of $\SL_2(\bbQ_v)$} \label{S1}

We will describe the double cover $G(\bbQ_v)$ in \eqref{eq1}.
If we fix a section $\bfs: \SL_2(\bbQ_v)\rightarrow G(\bbQ_v)$ then
$G(\bbQ_v)$ can be identified with the set $\SL_2(\bbQ_v) \times
\mu_2$, and the group law on $G(\bbQ_v)$ is given by
\[
(g_1, \epsilon_1) (g_2, \epsilon_2) = (g_1 g_2, \epsilon_1 \epsilon_2
\sigma_v(g_1,g_2))
\]
where $\sigma_v(g_1,g_2)$ is a cocycle which depends on
$\bfs$. Following  \cite{G2} we make
the following choice of the cocycle $\sigma_v$.  Let $(\, , \, )_v$ be
the Hilbert symbol over $\bbQ_v$. For $g = \begin{pmatrix} a & b \\ c
  & d
\end{pmatrix} \in \SL_2(\bbQ_v)$ we define
\begin{eqnarray*}
x(g) = \left\{
\begin{array}{ll}
c & \mbox{if } c \neq 0, \\
d & \mbox{if } c = 0;
\end{array}
\right. & \mbox{ and } &
s(g) = 
\left\{ 
\begin{array}{ll}
(c,d)_v & \mbox{if $v$ is a finite prime, $cd \neq 0$ }  \\
& \mbox{and ord$(c)$ is odd,} \\ 1 & \mbox{otherwise} .
\end{array}
\right.
\end{eqnarray*}
Then 
\[
\sigma_v(g_1,g_2) = (x(g_1g_2)x(g_1), x(g_1g_2) x(g_2))_v s(g_1)
s(g_2) s(g_1 g_2).
\]
                                   
An advantage of this particular section is that 
$K_p = \bfs(\SL_2(\bbZ_p))$ is a subgroup 
in $G(\bbQ_p)$ if $p\neq 2$. If $p = 2$, we define
\[
K_1(4) = \left\{ (
\begin{pmatrix}
a & b \\
c & d
\end{pmatrix}, 
1) \in \SL_2(\bbZ_2) \times \{ \pm 1 \} : a \in 1 + 4 \bbZ_2, c \in 4
\bbZ_2 \right\}.
\]
By Proposition 2.14 in \cite{G2}, $K_1(4)$ is a compact subgroup of
$G(\bbQ_2)$.

A smooth representation of $G(\bbQ_v)$ is call {\it genuine} if
$\mu_2$ acts non-trivially. If $p$ is an odd prime number, a smooth
genuine representation of $G(\bbQ_p)$ is called {\it unramified} if it
contains a vector fixed by $K_p$. A vector fixed by $K_p$ is called a
{\it spherical} vector. 

If $p = 2$, a smooth genuine representation is called {\it tamely
  ramified} if it contains a vector fixed by $K_1(4)$. Unfortunately
$\SL_2(\bbZ_2)$ does not split in $G(\bbQ_2)$ so we could not define
spherical vectors in the same manner as those for odd primes. The
objective of this paper is to motivate and define spherical vectors of
genuine representations of $G(\bbQ_2)$. 


\smallskip

We set up some notations for the later sections. For $u \in \bbQ_v$
and $t \in \bbQ_v^\times$, we define the following elements in
$\SL_2(\bbQ_v)$:
\[
\underline{x}(u) = \begin{pmatrix}
1 & u \\
0 & 1
\end{pmatrix}, \, 
\underline{y}(u)=\begin{pmatrix}  
1 & 0 \\
u & 1
\end{pmatrix}, \, \underline{w}(t) = \begin{pmatrix}
0 & t \\
-t^{-1} & 0
\end{pmatrix} \mbox{ and }
\underline{h}(t) = \begin{pmatrix}
t & 0 \\
0 & t^{-1}
\end{pmatrix}.
\] 
Let $x(u) = \bfs(\underline{x}(u))$, $y(u) =
\bfs(\underline{y}(u))$, $w(t) = \bfs(\underline{w}(t))$ and $h(t) =
\bfs(\underline{h}(t))$ in $G(\bbQ_v)$. Note that 
\[
h(t) h(s)=h(ts) (t,s)_v. 
\] 
Let $N = \{ x(u) : u \in \bbQ_v\}$, $\bar{N} = \{ y(u) : u\in \bbQ_v
\}$ and $T$ be the subgroup of $G$ generated by elements $h(t)$.


\section{Hecke algebra at $p = 2$} \label{S3}

We fix $p = 2$ throughout Sections \ref{S3} to \ref{S6}. We will
denote $G(\bbQ_2)$ by $G$ and $K_1(4)$ by $K_1$.  The objective of
these sections is to classify genuine
representations of $G$ containing a non-zero vector fixed by $K_{1}$.
  
Let $M$ be the center of $G$. It is a cyclic group of order $4$
generated by $h(-1)$.  (Note that $h(-1)h(-1) = (-1,-1)_2 = -1 \in
\mu_2$.) Thus, a genuine central character $\gamma$ is determined by
its value on $h(-1)$, which is a fourth root of 1. Let $K$ and $K_0$
be the open compact subgroups in $G$ equal to the inverse images of
$\SL_2(\bbZ_2)$ and
\[
\left\{ 
\begin{pmatrix}
a & b \\
c & d
\end{pmatrix}
\in \SL_2(\bbZ_2)  : c \in 4 \bbZ_2
 \right\}
\]
respectively. Let $K(4) \subset K_1$ denote the principal congruence
subgroup.  It is the image under the section $\bfs$ of the subgroup of
$\SL_2(\bbZ_2)$ consisting of matrices congruent to $1$ modulo~$4$.
We have $K \supset K_0 \supset K_1 \supset K(4)$ and $K_0=M\times
K_1$. We extend the central character $\gamma$ to $K_0$, so that it is
trivial on $K_1$.  Given a smooth representation $(\pi,V)$ of $G$, we
denote
\[ 
V^\gamma := \{ v \in V : \pi(k_0) v = \gamma(k_0) v \mbox{ for all }
k_0 \in K_0 \}.
\] 
Let $\calR(G,\gamma)$ denote the category of admissible smooth
(necessarily genuine) representations $V$ of $G$ such that $V^\gamma$
generates $V$ as a $G$-module.

Next we define the corresponding Hecke algebra.  Let $C_c(G)$ denote
the set of locally constant, compactly supported functions on $G$. Let
\[
H(\gamma) = \{ f : C_c(G) : f(k_0 g k_0') = \overgamma(k_0) f(g)
\overgamma(k_0') \mbox{ for all } k_0, k_0' \in K_0 \}.
\]
For $f_1, f_2 \in H(\gamma)$, we define 
\[
f_1 \cdot f_2(g_0) = \int_G f_1(g) f_2(g^{-1} g_0) dg = \int_G f_1(g_0 g)
f_2(g^{-1}) dg
\]
where $dg$ is the Haar measure on $G$ such that the measure of $K_0$
is 1. Then $H(\gamma)$ is a $\bbC$-algebra.  For $f \in
H(\gamma)$ and $v \in V$, we have
\[
\pi(f) v= \int_G f(g) \pi(g) v dg \in V^\gamma.
\]
In this way $V^\gamma$ is a left $H(\gamma)$-module. Let
$\calR(H(\gamma))$ denote the category of finite dimensional left
$H(\gamma)$-modules. We have a functor $A : \calR(G,\gamma)
\rightarrow \calR(H(\gamma))$ given by $V \mapsto V^\gamma$.  Since
the group $K_{0}$ has a triangular decomposition
\[
K_0 = (K_0 \cap \bar{N}) (K_0 \cap T) (K_0 \cap N)
\]
the functor $A$ is an equivalence of categories. This follows, in
essence, from \cite{Ca}, Corollary~3.3.6 (see also \cite{Bo} and
Theorem 4.2 in \cite{BZ}).

\smallskip

Our immediate goal is to understand the structure of $H(\gamma)$. The
character $\gamma$ of the center $M$ extends to a character $\gamma$
of $T$ which is trivial on $K_1 \cap T$ and $\gamma(h(2^n)) = 1$ for
all $n \in \bbZ$. Let us abbreviate 
\[
\gamma(t) = \gamma(h(t)).
\]
We define
\[ 
\zeta = \frac{1+\gamma(-1)}{\sqrt{2}}.
\] 
Note that $\zeta$ is a primitive 8-th root of 1.  The 
character $\gamma$ of $T$ is invariant under conjugation
by $w = w(1)$.  We can now extend the character $\gamma$ from $T$ to
the normalizer $N_G(T)$ by defining $\gamma(w) = \zeta$.

We define some functions in $H(\gamma)$.  For $g$ in $N_G(T)$ we set
$X_g$ to be the function supported on $K_0 g K_0$ such that
\[
X_g(k_0 g k_0') = \overgamma(k_0) \overgamma(g) \overgamma(k_0')
\] 
for all $k_0, k_0' \in K_0$. Note that this definition depends only on
the image of $g$ in the affine Weyl group $W_a := N_G(T)/(T \cap
K_0)$.

\begin{prop} \label{P1}
 Functions ${X}_g$ for $g$ in $W_a$ form a basis of $H(\gamma)$. 
\end{prop}  
\begin{proof} 
We need first to determine the $K_0$-double cosets in $G$. This can be
easily determined in $\SL_2(\bbQ_2)$ using the row-column
reduction. In addition to $h(2^n)$ and $w(2^{-n})$ the double coset
representatives are:
\[ 
y(2), h(2^n) y(2), y(2)h(2^{-n}), y(2) w(2^{-n}), w(2^{-n})y(2)
\text{ and } y(2)w(2^{-n})y(2)
\] 
where $n\geq 1$ in all cases. We claim that the Hecke algebra is not
supported on these cosets. 

\begin{lemma} 
The commutator of $x(2)$ and $y(2)$ modulo the principal congruence
subgroup $K(4)$ is equal to $-1\in \mu_2$.
\end{lemma}

\begin{proof} 
This can be easily checked using the multiplication rule. It also
follows from applying Corollary 2.9 in \cite{St} to the ring $A = \bbZ/4\bbZ$, 
 \end{proof} 

Now we can easily finish the proof of proposition. Indeed if $f$ is in
$H(\gamma)$ then
\[
f(y(2))=f(y(2)x(2))=-f(x(2)y(2))=-f(y(2))
\]
by the above lemma. This implies that $f$ must vanish on $y(2)$. Other
cases are dealt with in the same manner.
\end{proof}

Let $\ell : N_G(T)\rightarrow \bbZ$ be defined by $\ell(g)= \log_2(n)$
where $n$ is the number of left (or right) $K_0$-cosets in the double
coset $K_0 g K_0$. In other words, the volume of $K_0gK_0$ is
$2^{\ell(g)}$. For example, $w(2^{-1})$ normalizes $K_0$, so
$\ell(w(2^{-1}))=1$.

\begin{prop} \label{P3}
For every integer $n$ we have $\ell(h(2^n))= 2|n|$ and
$\ell(w(2^{-n}))=2|1-n|$. More precisely, we have the following
decompositions of double co-sets:
\begin{enumerate}[(i)]
\item If $n \geq 0$,
\[ 
K_0 h(2^n) K_0 = \bigcup_{u \in \bbZ/2^{2n} \bbZ} x(u) h(2^n) K_0 =
\bigcup_{u \in \bbZ/2^{2n} \bbZ} K_0 h(2^n) y(4u).
\]
\item If $n \geq 1$,
\[ 
K_0 h(2^{-n}) K_0 = \bigcup_{u \in \bbZ/2^{2n} \bbZ} y(4u) h(2^{-n})
K_0 = \bigcup_{u \in \bbZ/2^{2n} \bbZ} K_0 h(2^{-n}) x(u).
\]
\item If $n\geq 0$,
\[ 
K_0 w(2^n) K_0 = \bigcup_{u \in \bbZ/2^{2n+2}\bbZ} x(u) w(2^n) K_0 =
\bigcup_{u \in \bbZ/2^{2n+2}\bbZ} K_0 w(2^n) x(u).
\] 
\item If $n\geq 1$,
\[
K_0 w(2^{-n}) K_0 = \bigcup_{u \in \bbZ/2^{2n-2} \bbZ} y(4u) w(2^{-n})
K_0 = \bigcup_{u \in \bbZ/2^{2n-2} \bbZ} K_0 w(2^{-n}) y(4u).
\] 
\end{enumerate}
\end{prop}    

\begin{proof} This  is an easy consequence of the decomposition 
$K_0=(K_0\cap \bar{N})(K_0\cap T)(K_0\cap N)$. Details are left to the reader. 
\end{proof} 
    
We record the following tautological lemma: 

\begin{lemma}\label{L4}
 Let $g_1$ and $g_2$ be two elements in $N_G(T)$. If $\ell(g_1 g_2) = 
\ell(g_1) + \ell(g_2)$ then $X_{g_1} \cdot X_{g_2}=X_{g_1g_2}$. \qed
\end{lemma}  

\smallskip 
Let 
\[ 
\begin{cases} T_n = X_{h(2^n)}\\
U_n = X_{w(2^{-n})}. 
\end{cases} 
\] 

\begin{prop} \label{P5}  
Let $T_w = \sqrt{2}^{-1} U_0$. We have the following identities where
$m, n$ are any integers unless specified otherwise.
\begin{enumerate}[(i)]
\item $(T_w+1)(T_w-2) = 0$.

\item $U_1 \cdot  U_1 = 1$.
  
\item If $m, n \geq 0$, or $m, n \leq 0$ then $T_m \cdot T_n = T_{m+n}$.

\item $U_1\cdot  T_n = U_{n+1}$ and $T_n \cdot  U_1 = U_{1-n}$.

\item $U_1 \cdot U_n = T_{n-1}$ and  $U_n \cdot U_1 = T_{1-n}$.
\end{enumerate}
\end{prop}

\begin{proof} 
All statements except the first follow from Lemma \ref{L4}.
For (i) we need to show $T_w^2 = T_w\cdot T_w = T_w + 2$.  Since
$T_w^2$ is supported in $K$ this is equivalent to $T_w^2(1) = 2$ and
$T_w^2(w(1)) = T_{w}(w(1))$.  Suppose $f_1, f_2 \in H(\gamma)$ where
$f_1$ is supported on $K_0 r K_0 = \bigcup_{i=1}^s r_i K_0$
(disjointed union). Then
\begin{eqnarray*}
f_1 \cdot f_2(g) & = &  \sum_{i = 1}^s
f_1(r_i) f_2(r_i^{-1}g).
\end{eqnarray*}
We can apply this observation to $f_{1}=f_{2}=T_{w}$.  Proposition
\ref{P3} (the case $n=0$ in (iii)) gives a decomposition of
$K_{0}w(1)K_{0}$ into single single cosets. Hence
\[
T_w^2(g) = \sum_{u \, {\mathrm{mod}} \, 4} T_{w}(x(u)w(1))\cdot 
T_w(w(-1) x(-u) g). 
\]
If $g =1$, this gives $T_{w}^{2}(1)=4T_{w}(w(1))\cdot T_{w}(w(-1))$.
Since $T_{w}(w(1))=2^{-1/2} \overzeta$ and $T_{w}(w(-1))=2^{-1/2}
\zeta$, we obtain that $T_{w}^{2}(1)=2$.  If $g = w(1)$, then
\[
T_w^2(w(1)) = T_{w}(w(1))
\sum_{u \, {\mathrm{mod}} \, 4} T_w(y(u)).
\] 
If $u=0$ or $2$ then $y(u)$ is not in $K_{0}w(1)K_{0}$ and $T_w(y(u))
= 0$. If $u = \pm 1$, then $y(u) = x(u) w(-u) x(u)$ and we can rewrite
\[
T_w^2(w(1)) = T_{w}(w(1))[T_w(w(1)) + T_{w}(w(-1))]=T_{w}(w(1)).
\] 
This proves (i).
\end{proof}

Here is the main  result of this section. 

\begin{thm} \label{T6}
The Hecke algebra $H(\gamma)$ is generated by $T_w$ and $U_1$ as an
abstract $\bbC$-algebra modulo the relations
\begin{enumerate}[(a)]
\item $(T_w-2)(T_w+1) = 0$ and
\item $U_1^2 = 1$.
\end{enumerate}
\end{thm}

\begin{proof} 
  Suppose $H$ is the abstract algebra generated by $U_{0}=\sqrt{2}T_w$
  and $U_{1}$ modulo the relations (a) and (b). We have a natural
  homomorphism of $\bbC$-algebras $B: H\rightarrow H(\gamma)$.  By
  Proposition \ref{P1}, $H(\gamma)$ is spanned by $T_n$ and $U_n$ and
  by Proposition \ref{P5}, these elements are generated by $U_0$ and
  $U_1$. This shows that $B$ is surjective. It remains to show that
  $B$ is injective.  Suppose $h \in H$ is in the kernel of $B$. Since
  $U_{0}$ and $U_{1}$ satisfy quadratic relations, $h = \sum_i c_i
  u_i$ where $c_i \in \bbC$ and $u_i \in H$ is of the form $U_{1} U_0
  U_{1} U_0 \ldots$ or $U_0 U_{1} U_0 U_{1} \ldots$.  Since
  $U_{0}U_{1} = T_{1}$, $B(u_i)$ is either $T_n$, $T_{n}U_{1} =
  U_{1-n}$, $U_{1}T_n = U_{n+1}$, or $U_{1}T_n U_{1} = T_{-n}$. These
  elements have disjointed supports as functions in $H(\gamma)$.
  Therefore $B(h) = \sum_i c_i B(u_i) = 0$ implies that $c_i = 0$ and
  $h = 0$. This proves that $B$ is an injection and Theorem \ref{T6}.
\end{proof} 

\smallskip

We now give two consequences of Theorem \ref{T6}: 

\begin{prop} \label{P7}
The element $Z := \frac{T_1}{2} + (\frac{T_1}{2})^{-1}$ belongs to the
center of $H(\gamma)$.
\end{prop}

\begin{proof}
By Proposition \ref{P5}, $T_1$ and $U_{1}$ generate
$H(\gamma)$. Clearly $Z$ commutes with $T_1$. It suffices to show that
$Z$ commutes with $U_1$. Since $T_1 = U_0 U_{1}$ we can use quadratic
relations satisfied by $U_{0}$ and $U_1$ to write
\begin{equation} \label{eq2}
2Z = U_0 U_{1} + U_1U_{0} - 2^{1/2} U_{1}.
\end{equation}
Hence $Z$ commutes with $U_{1}$. This proves the proposition.
\end{proof}

\begin{prop} \label{P8}
For $n \geq 0$, $T_n$ is an invertible element in the algebra
$H(\gamma)$.
\end{prop}

\begin{proof} 
Note that the quadratic relations satisfied by $U_0$ and $U_1$ imply
that $U_0$ and $U_1$ are invertible. Since $T_1=U_0 U_1$, $T_1$ is
also invertible. Hence $T_n = T_1^n$ is invertible.
\end{proof}

Suppose $(\pi,V)$ is a representation in $\calR(G,\gamma)$.  Let
$(V_N)^\gamma = \{ v \in V_N : \pi_{V_N}(t) v = \gamma(t) v \mbox{ for
  all } t \in K_0 \cap T \}$.  The invertibility of $T_{n}$ implies
(see Lemma 4.7 in \cite{Bo}):

\begin{cor} \label{C9}
 Suppose $(\pi,V)$ is a representation in $\calR(G,\gamma)$.  Then the
 canonical map $V^\gamma \rightarrow (V_N)^\gamma$ is a bijection. In
 particular $V_N$ is nonzero, and $V$ cannot be a supercuspidal
 representation. \qed
\end{cor}

\section{Local Shimura correspondence} \label{S4}

Let $G' = \PGL_2(\bbQ_2)$. Let $I$ be its Iwahori subgroup and let
$H'$ be its Iwarhori-Hecke algebra. Let  $T_w'$ and $U_{1}'$ denote
the characteristic functions of
\[
I
\begin{pmatrix} 
0 & 1 \\ 
1 & 0 
\end{pmatrix} 
I \text{ and } I
\begin{pmatrix} 
0 & 1 \\
p & 0 
\end{pmatrix} I
\]
respectively. Then $H'$ is the abstract $\bbC$-algebra generated by
$T_w'$ and $U'_{1}$ satisfying the same relations as Theorem~\ref{T6}
(a) and (b) (see \cite{Ma}). This gives the next corollary.

\begin{cor} \label{C10}
  The Hecke algebras $H(\gamma)$ and $H'$ are isomorphic
  $\bbC$-algebras. \qed
\end{cor}

Let $\calR(H')$ denote the category of finite dimensional
representations of $H'$. Let $\calR(G',I)$ denote the category of
admissible smooth representations $V$ of $G'$ such that $V^{I}$
generates $V$ as a $G'$-module.  By \cite{Bo} and \cite{BZ}, the
functor $V \mapsto V^{I}$ is an equivalence of categories from
$\calR(G',I)$ to $\calR(H)$.  The isomorphism in Corollary \ref{C10}
establishes an equivalence of categories between $\calR(H(\gamma))$
and $\calR(H')$. Hence the following four categories are equivalent:
\[
\calR(G,\gamma) \simeq \calR(H(\gamma)) \simeq \calR(H')
\simeq \calR(G',I). 
\]
If $V$ is a representation in $\calR(G,\gamma)$, then we call the
corresponding representation in $\calR(G',I)$ the {\it local Shimura
  lift} of $V$. We denote it by $\rSh(V)$.

\begin{prop} \label{P11}
  Let $V$ be a representation in $\calR(G,\gamma)$.  Then the
  following statements are equivalent.

\begin{enumerate}[(i)]
\item The local Shimura lift $\rSh(V)$ is a spherical representation
  of $G'$.

\item 
The action of $T_w'$ on $\rSh(V)^{I}$ has an eigenvalue $2$.

\item 
The action of $T_w$ on $V^\gamma$ has an eigenvalue $2$. 
\end{enumerate}
\end{prop}

\begin{proof}
  The projection map to $G'(\bbZ_2)$-fixed vectors in $\rSh(V)$ is
  given by $\frac{1}{3}(T_w' + 1)$, since $T_w'+1$ is the
  characteristic function of $G'(\bbZ_2)$ and the volume of
  $G'(\bbZ_2)$ is 3. It follows that a $G'(\bbZ_2)$-fixed vector is
  an eigenvector of $T_w'$ with eigenvalue $2$.  This proves the
  equivalence of (i) and (ii). The equivalence of (ii) and (iii)
  follows from Corollary~\ref{C10}.
\end{proof}

The above proposition motivates the following definition.

\subsection*{Definition}
Let $V$ be a smooth representation of $G$. An eigenvector of $T_w$ in
$V^\gamma$ with an eigenvalue 2 is called a {\it $\gamma$-spherical}
vector. The representation is called a {\it $\gamma$-unramified} or
{\it $\gamma$-spherical} representation if it contains a
$\gamma$-spherical vector.

\smallskip

\section{Pseudo-spherical representation of $K$ at $p = 2$} \label{S5}

We retain the notations in Sections \ref{S3} and \ref{S4} where $p =
2$.  In the previous section we defined a representation $V$ of $G$ to
be unramified if $V^{\gamma}\neq 0$ and $T_w$ has an eigenvalue $2$ on
$V^{\gamma}$. In this section we shall reinterpret this condition in
terms of representations of $K$.  We shall see that $K$ has only two
irreducible representations $E$ such that $E^{\gamma}\neq 0$. For both
representations $E^{\gamma}$ is one dimensional and they are
distinguished by the action of $T_w$ on $E^{\gamma}$.  That eigenvalue
can be either $2$ or $-1$, so we shall use the eigenvalue to denote
the representations by $V(2)$ and $V(-1)$.  Their dimensions are $2$
and $4$, respectively. Thus, a representation of $G$ is unramified if
and only if it contains the two-dimensional $K$-type $V(2)$, which we
may call a pseudo-spherical type.

If $E^{\gamma}\neq 0$ then, by Frobenius reciprocity, the $K$-type $E$
is a summand of a six dimensional induced representation
\[
I_K(\gamma) := \Ind_{K_0}^K \gamma = \{\phi : K \rightarrow \bbC :
\phi(k_0 k) = \gamma(k_0) \phi(k) \mbox{ for all } k \in K, k_0 \in
K_0 \}.  
\]  
Here the group $K$ acts on it by right translation. We denote this
action by $\pi_R$. Let $H_K(\gamma)$ denote the subalgebra of
$H(\gamma)$ consisting of functions supported on $K$. We have the
action of $H(\gamma)$ on $I_{K}(\gamma)^{\gamma}$, also denoted by
$\pi_{R}$.  By Proposition \ref{P1}, $H_K(\gamma) = \bbC 1 \oplus \bbC
T_w$ and it is a commutative subalgebra.  The algebra $H_{K}(\gamma)$
is anti-isomorphic to the algebra $H_{K}(\bar{\gamma})$ via the map
$f\mapsto \hat{f}$ where
\[
\hat{f}(g) = f(g^{-1}).
\]

For $f \in H_K(\bar{\gamma})$ and $\phi \in I_K(\gamma)$, we set
\[
(\pi_L(f) \phi)(g) = \int_K f(k) \phi(k^{-1} g) d k \ \ 
\mbox{for all } g \in K. 
\]
This action commutes with the right action $\pi_R$ of $K$ on
$I_K(\gamma)$ and
\[ 
H_K(\bar\gamma) = \End_K(I_K(\gamma)).
\] 
Note that $I_{K}(\gamma)^{\gamma}=H(\bar\gamma)$. The actions
$\pi_{L}$ and $\pi_{R}$ of $H(\bar{\gamma})$ and $H(\gamma)$ on
$I_{K}(\gamma)^{\gamma} = H(\bar\gamma)$ are related by
$\pi_{L}(\hat{f})=\pi_{R}(f)$.

We define the functions $\phione := \frac{1}{3}(2 - T_w)$ and
$\phitwo := \frac{1}{3}(T_w+1)$ in $H_K(\gamma)$. Then $\{
\phione, \phitwo \}$ is a basis of idempotents of $H_K(\gamma)$.

For $j = -1,2$, let $V(j) = \pi_L(\hat \phij) I_K(\gamma)$. In other
words $V(j)$ is the eigenspace of $\pi_L(\hat{T}_w)$ on $I_K(\gamma)$
corresponding to the eigenvalue~$j$. Note that $\hat{F}_j \in V(j)$
and $\pi_R(T_w) \hat{F}_j = j \hat{F}_j$. In particular $\hat{F}_2$ is
a $\gamma$-spherical vector.

\begin{prop} \label{P12}
\begin{enumerate}[(i)]
\item 
We have $I_K(\gamma) = V(-1) \oplus V(2)$ where each summand is an
irreducible representation of $K$. 

\item 
We have $\dim V(-1) = 4$ and $\dim V(2) = 2$.

\item The $K$-submodule $V(2)$ contains a $\gamma$-spherical vector
  $\hat{F}_2$. The space of $\gamma$-spherical vectors is one
  dimensional.

\item 
The $K$-submodule $V(-1)$ does not have any $\gamma$-spherical vector.
\end{enumerate}
\end{prop}

\begin{proof}
Since $\dim \End(I_K(\gamma)) = 2$, both $V(-1)$ and $V(2)$ are
irreducible $K$-modules. This proves (i).

In order to compute the dimensions of $V(-1)$ and $V(2)$ we need the
following lemma.

\begin{lemma}
The operator $\pi_L(\hat{T}_w)$ as an element in $\End_K(I_K(\gamma))$
has trace $0$.
\end{lemma}

\begin{proof}
  For $g \in K$, let $\phi_g \in I_K(\gamma)$ such that $\phi_g$ is
  supported on $K_0g$ and $\phi_g(k_0g) = \gamma(k_0)$.  Let $S$ be a
  set of representatives of $K_0 \backslash K$, then $\{ \phi_g :
  g \in S\}$ is a basis of $I_K(\gamma)$.  In order to prove the
  lemma, it suffices to show that $(\pi_L(\hat{T}_w) \phi_g)(g) = 0$.
  Indeed, this shows that the matrix of $\pi_{L}(\hat{T}_{w})$ in the
  basis $\phi_{g}$ has all diagonal entries equal 0.  Note that
  $\pi_L(T_w) \phi_g$ is supported on $K_0 w(1) K_0 g$. If
  $(\pi_L(\hat{T}_w) \phi_g)(g) \neq 0$, then $g \in K_0 w(1) K_0 g$
  and $1 \in K_0 w(1) K_0$. Since $K_{0}\neq K_{0} w(1) K_{0}$, this is
  a contradiction. The lemma is proved.
\end{proof}

We have $\dim V(2) + \dim V(-1) = \dim I_K(\gamma) = [K:K_0] = 6$. By
the above lemma, $2\dim V(2) - \dim V(-1) = 0$. This implies $\dim
V(-1) = 4$ and $\dim V(2) = 2$ and proves Proposition \ref{P12}(ii).
We have $I_K(\gamma)^\gamma = H_K(\bar\gamma)$ and $\pi_R(F_j)
I_K(\gamma) = \bbC \hat{F}_j$ for $j = -1,2$.  The vector
$\hat{F}_2$ is $\gamma$-spherical while $\hat{F}_{-1}$ is
not. This proves Parts (iii) and (iv).
\end{proof}

\begin{thm}
  A smooth representation $V$ of $G$ with central character $\gamma$
  is {\it $\gamma$-unramified} if and only if there is a nontrivial
  $K$-module homomorphism $l : V(2) \rightarrow V$.  A vector in $V$
  which is a scalar multiple of $l(\hat{F}_2)$ is a
  $\gamma$-spherical vector of $V$.
\end{thm}

\begin{proof}
A $\gamma$-spherical vector in $V$ would generated a representation of
$K$ where every irreducible $K$-submodule is isomorphic to an
irreducible submodule of $I_K(\gamma)$. Now the theorem follows from
Proposition \ref{P12}.
\end{proof}

\section{Unramified principal series representations at $p = 2$} 
\label{S6}

In this section, we continue to assume $p = 2$ and notations as in
Sections \ref{S3} to~\ref{S5}.  We will show that $\gamma$-unramified
representations appear as submodules of principal series
representations.

\smallskip

We recall the character $\gamma$ of $T$ in Section \ref{S3}. Let
$(\pi_s, I(\gamma,s))$ be the normalized induced principal series
representation where $I(\gamma,s)$ is the set of smooth functions
$\phi : G \rightarrow \bbC$ satisfying
\[
\phi(\epsilon x(u)h(t)  g) = \epsilon \gamma(t) |t|^{s+1} \phi(g)
\]
for all $\epsilon \in \mu_2$, $u \in \bbQ_{2}$ and $t \in \bbQ_2^\times$. 
The group $G$ acts by left translation $(\pi_s(g) \phi)(g') =
\phi(g'g)$.
 
\begin{prop} \label{P15}
An irreducible $\gamma$-unramified representation $V$ is isomorphic to
a submodule of some $I(\gamma,s)$.
\end{prop}

\begin{proof}
By Corollary \ref{C9}, $(V_N)^\gamma$ is nonzero. Hence there is a
nontrivial $T$-homomorphism $V_N \rightarrow \gamma \nu^{s+1}$ for
some $s \in \bbC$. Here $\nu$ is the character $\nu(\underline{h}(t))
= |t|$.  By Frobenius reciprocity, there is a nontrivial map $V
\rightarrow I(\gamma,s)$ which is an injection because $V$ is
irreducible.
\end{proof}

We recall that $K(4)$ is the principal congruence subgroup in
$K_{1}$. Restricting functions $\phi$ in $I(\gamma,s)$ to $K$ gives a
natural isomorphism of $K$-modules
\[
l: I_{K}(\gamma) \rightarrow I(\gamma,s)^{K(4)}.
\]

\begin{thm} 
The $K$-types $V(2)$ and $V(-1)$ appear with multiplicity one in
$I(\gamma,s)$. The space $I(\gamma,s)^{\gamma}$ is $2$-dimensional. It
is spanned by $l(\hat{F}_{2})$ and $l(\hat{F}_{-1})$. \qed
\end{thm}

We will describe a scalar multiple $\phi_{j}$ of $l(\hat{F}_{j}) \in V_j$ which
is more convenient for later calculations. Let $d_2 = 1$ and $d_{-1} =
-2$, and define $\phi_j$ be the unique vector in $I(\gamma,s)$ whose
restriction to $K$ is given by
\[
\phi_j(k) = \left\{ 
\begin{array}{ll} 
d_j \gamma(k) & \mbox{ if } k \in K_0, \\ 2^{-1/2} \zeta
\gamma(k_0 k_0') & \mbox{ if } k = k_0 w(1) k_0' \in K_0 w(1) K_0, \\ 0 &
\mbox{ otherwise.}
\end{array}
\right.
\]

\smallskip

We define an intertwining map $M(s) : I(\gamma,s) \rightarrow
I(\gamma,-s)$ by
\[
(M(s) \phi)(g) = \int_{\bbQ_2} \phi(w(1) x(u) g) du  
\]
where $g$ is in $G$ and $du$ is the Haar measure on $\bbQ_2$ such that
the measure of $\bbZ_2$ is 1.

\begin{prop} \label{P17}
We have
\[
M(s) \phi_2 = \frac{\zeta}{\sqrt{2}} \left(\frac{1-\frac{1}{2}
  (2^{-2s} )}{1-2^{-2s} } \right) \phi_2 \mbox{ \ and \ } M(s) \phi_{-1} =
-\frac{\zeta}{2 \sqrt{2}} \left(\frac{1-2 (2^{-2s} )}{1-2^{-2s} } \right)
\phi_{-1}.
\]
\end{prop}

\begin{proof}
Since the vector $\phi_j$ is unique up to a scalar in $I_K(\gamma)$,
$M(s) \phi_j = c \phi_j$ for some $c \in \bbC$. It remains to
determine $c = d_j^{-1} M(s) \phi_j(1)$.

If $u \not \in \bbZ_2$, then $w(1) x(u) = (-1,u)_{2}\cdot x(-u^{-1})
h(u^{-1}) y(u^{-1})$. We write $u^{-1} = 2^{m} v$ where $v \in
\bbZ_2^\times$ and $m\geq 1$.  Recall that $\gamma(t) =
\gamma(h(t))$. Then
\begin{eqnarray*} M(s) \phi_j(1) & = &
    \int_{\bbZ_2} \phi_j(w(1)x(u)) du + \sum_{m = 1}^\infty 2^{m-1}
    \int_{\bbZ_2^\times} (-1,2^m v)_2 \phi_j(h(2^m v) y(2^m
    v)) d^\times v \\ & = & 2^{-1/2} \zeta +
    \sum_{m = 1}^\infty 2^{-ms-1} \int_{\bbZ_2^\times} (-1,v)_2
    \gamma(2^m v) \phi_j(y(2^m v)) d^\times v
\end{eqnarray*}
where $d^\times v$ is the Haar measure of $\bbZ_2^\times$ with total
measure 1.  Now $\phi_j(y(2^m v)) = 0$ if $m = 1$ and it is equal to
$1$ if $m \geq 2$. Since $\gamma(2^m v) = \gamma(2^{m})
\gamma(v)(2^{m},v)_{2}$ and $\gamma(2^{m}) = 1$, we can rewrite
\begin{eqnarray*}
M(s) \phi_j(1) & = & 2^{-1/2} \zeta + d_j \sum_{m = 2}^\infty
2^{-ms-1} \int_{\bbZ_2^\times} (2,v)_2^m (-1,v)_2
\gamma(v) d^\times v \\ & = & 2^{-1/2} \zeta + d_j
\sum_{m = 2}^\infty 2^{-ms-1} \frac{1}{4} \sum_{v \in
  (\bbZ/8\bbZ)^\times} (2,v)_2^m (-1,v)_2 \gamma(v).
\end{eqnarray*}
The sum $\sum_{v \in (\bbZ/8\bbZ)^\times}$ on the right is zero if
$m$ is odd, and equals $\sqrt{2} \zeta$ if $m$ is even. Finally adding
up all the terms gives the constant $c$ and the lemma.
\end{proof}

\smallskip

Let $s_0 = \frac{1}{2}$ or $\frac{1}{2} + \frac{i \pi}{\log 2}$. From
the above proposition, $\phi_{-1}$ lies in the kernel of $M(s_0)$ so
$I(\gamma, s_0)$ is reducible.  Indeed $I(\gamma, s_0)$ has a unique
irreducible quotient which is an even Weil representation.

\subsection*{Definition}
Let $s_0 = \frac{1}{2}$ or $\frac{1}{2} + \frac{i \pi}{\log 2}$. The
kernel of $M(s_0)$ is called the {\it Steinberg} representation of
$G(\bbQ_2)$.  We shall denote this representation by $St(\epsilon)$
where $\epsilon=\pm 1$ such that $2^{s_0}=\epsilon \sqrt{2}$. 

\smallskip

We claim that $St(\epsilon)$ is an irreducible representation of
$G(\bbQ_2)$. Indeed by Section 6 in \cite{LS}, for every $s \in \bbC$,
we have
\begin{equation} \label{eq3}
I(\gamma,s)^{\rmss}_N \cong \gamma|\cdot|^{s+1} \oplus
\gamma|\cdot|^{-s+1} 
\end{equation}
where $I(\gamma,s)^{\rmss}_{N}$ is the semi-simplification of
$I(\gamma,s)_{N}$ as a $T$-module. Hence $I(\gamma,s)$ has at most
length 2. The claim now follows because $St(\epsilon)$ is a proper
submodule of $I(\gamma, s_0)$. Also see Section 7 of \cite{Sa}.

 
\begin{cor} \label{C18}
The even Weil representation contains the irreducible $K$-module
$V(2)$. It is a $\gamma$-unramified representation. The Steinberg
representation contains the irreducible $K$-module $V(-1)$. \qed
\end{cor}

\begin{prop} \label{P19}
 Let $Z = \frac{T_{1}}{2}+(\frac{T_{1}}{2})^{-1}$ be the central
 element in the Hecke algebra $H(\gamma)$ as in Proposition
 \ref{P7}. Then $\pi_{s}(Z)$ acts on $I(\gamma,s)^\gamma$ as the
 scalar $2^{s}+2^{-s}$.
\end{prop}

\begin{proof} 
By Corollary \ref{C9}, the natural projection of $I(\gamma,s)$ on
$I(\gamma,s)_N$ gives an isomorphism of $I(\gamma,s)^{\gamma}$ and
$I(\gamma,s)_N$.  From the decomposition of $K_0h(2)K_0$ into single
$K_0$-cosets (Proposition \ref{P3}(i)) it follows that the action of
$T_1$ on $I(\gamma,s)^{\gamma}$ corresponds to the action of $4\cdot
\pi_{s,N}(h(2))$ on $I(\gamma,s)_N$. By \eqref{eq3} the eigenvalues of
$\frac{T_1}{2}$ are $2^s$ and $2^{-s}$. This proves the proposition.
\end{proof} 

\begin{cor} \label{C20}
  An irreducible $\gamma$-unramified representation is uniquely
  determined by the eigenvalue of the action of $Z$ on its
  $\gamma$-spherical vector.
\end{cor}

\begin{proof}
  Suppose the irreducible $\gamma$-unramified representation is a
  subquotient of both $I(\gamma,s)$ and $I(\gamma,s')$. Then by
  Proposition \ref{P19}, $2^s + 2^{-s} = 2^{s'} + 2^{-s'}$ which
  implies $2^s = 2^{s'}$ or $2^s = 2^{-s'}$. By Proposition \ref{P17}
  both $I(\gamma,s)$ and $I(\gamma,-s)$ have the same irreducible
  $\gamma$-unramified subquotient. This proves the corollary.
\end{proof}

\begin{cor} \label{Steinberg} 
 The Steinberg representation $St(\epsilon)$ corresponds to the 
one dimensional representation of $H(\gamma)$ given by $T_w=-1$ and 
$U_1=-\epsilon$.
\end{cor} 
\begin{proof} We know that $T_w=-1$  
  on $St(\epsilon)^{\gamma}$. It remains to compute the action of
  $U_1$.  Since $St(\epsilon)$ is a subquotient of $I(\gamma, s_0)$
  where $2^{s_0}= \epsilon \sqrt{2}$, the central element $Z$ acts on
  $St(\epsilon)$ by the scalar $\epsilon(2^{1/2}+2^{-1/2})$. By
  \eqref{eq2} we have $2^{1/2}Z=T_w U_1+ U_1T_w - U_1$. Hence
  $U_1=-\epsilon$ as claimed.
\end{proof} 
\smallskip 

Let $V$ be an irreducible $\gamma$-unramified representation. By
Proposition \ref{P15}, we may assume that $V$ is the unique
$\gamma$-unramified subquotient of $I(\gamma,s)$ for some $s \in
\bbC$. By Proposition~\ref{P11}, its local Shimura lift $V' = \rSh(V)$
is an unramified irreducible representation of $G' =
\PGL_2(\bbQ_2)$. Let $B'$ be the Borel subgroup of $G'$. We may
realize $V'$ as the unramified irreducible subquotient of the
normalized induced principal series representation $(\pi_s', I'(t))$
with trivial central character. Here $I'(t) = \Ind_{B'}^{G'} \omega^t$
(normalized induction) where $\omega$ is the character
\[
\omega \begin{pmatrix} a_{1} & 0 \\ 0 & a_{2} \end{pmatrix} =
|a_{1}/a_{2}|.
\] 

\begin{thm} \label{T22}
If $V$ is the unique $\gamma$-unramified irreducible subquotient of
$I(\gamma,s)$, then its local Shimura lift $\rSh(V)$ is the unique
unramified irreducible subquotient of $I'(s)$.
\end{thm}

\begin{proof}
Assume that $\rSh(V)$ is a subquotient of $I'(t)$. By Proposition
\ref{P19} the central operator $Z$ in $H(\gamma)$ acts on
$I(\gamma,s)^\gamma$ by the scalar $2^s+2^{-s}$. The corresponding
operator $Z'$ in the algebra $H'$ acts on $I'(t)$ by $2^t+2^{-t}$.
Thus, $2^s + 2^{-s} = 2^t + 2^{-t}$. Solving the equation gives $2^s =
2^t$ or $2^s = 2^{-t}$. Both $I'(t)$ and $I'(-t)$ have the same
irreducible subquotients so we may set $s = t$.  This proves the
theorem.
\end{proof}

\begin{cor} \label{C23} 
The principal series representation $I(\gamma,s)$ is reducible if and
only if $s = \frac{1}{2}$ or $\frac{1}{2} + \frac{i \pi}{\log 2}$.
\end{cor}

\begin{proof}
Let $V$ be the $\gamma$-unramified irreducible subquotient of
$I(\gamma,s)$. Let $W$ be the unramified irreducible subquotient of
$I'(s)$. Then $V = I(\gamma,s)$ if and only if $\dim V^\gamma = 2$. By
Theorem \ref{T22}, $\dim V^\gamma = \dim W^I$. Now $\dim W^I = 2$ if
and only if $I'(s)$ is irreducible. Finally $I'(s)$ is irreducible if
and only if $s \neq \frac{1}{2}$ and $\frac{1}{2} + \frac{i \pi}{\log
  2}$.
\end{proof}

\section{Automorphic forms} \label{S7}

In this section we first review a connection between automorphic forms
and classical modular forms of half integral weight. This is mostly a
well known material that can be found in Chapters 2 and 3 of
\cite{G2}, and in \cite{Wa}. We then transfer the action of the Hecke
algebra $H(\gamma)$ to the setting of classical modular forms.

Let $\bbA = \prod_v \bbQ_v$ be the ring of adeles over $\bbQ$. We
recall $K_p$, $s(g)$ and the cocycle $\sigma_v$ defined in
Section~\ref{S1}.  Let $G(\bbA) = \SL_2(\bbA) \times \{ \pm 1 \}$ as a
set. For $g_1 = (g_{1,v})$, $g_2 = (g_{2,v}) \in \SL_2(\bbA)$ and
$\epsilon_1, \epsilon_2 \in \{ \pm 1\}$, the group law on $G(\bbA)$ is
given by
\[
(g_1, \epsilon_1) (g_2, \epsilon_2) = (g_1 g_2, \epsilon_1 \epsilon_2
\sigma(g_1,g_2))
\]
where $\sigma(g_1,g_2) = \prod_v \sigma_v(g_{1,v},g_{2,v})$. Then $\pr
: G(\bbA) \rightarrow \SL_2(\bbA)$ given by $\pr(g,\epsilon) = g$ is a
two-fold cover which splits over the subgroup $\SL_{2}(\bbQ)$.  Since
$\SL_{2}(\bbQ)$ is perfect this splitting is unique and it is given by
$\bfs_\bbQ: \SL_2(\bbQ) \rightarrow G(\bbA)$
\[
\bfs_\bbQ(g) = (g,s_\bbA(g))
\] 
where $s_\bbA(g) = \prod_v s(g_v)$.  

We also need a description of a maximal compact subgroup in $G(\bbR)$.
Let 
\[
\underline{k}(\theta) = \begin{pmatrix} \cos \theta & \sin \theta \\
  -\sin \theta & \cos \theta \end{pmatrix} \in \SL_2(\bbR)
\]
for $ -\pi <\theta \leq \pi$. Then $\underline{K}_\infty := \{
\underline{k}(\theta) : -\pi <\theta \leq \pi \}$ is a maximal compact
subgroup in $\SL_{2}(\bbR)$. Let $K_{\infty} = \{ k(\theta) : -2 \pi <
\theta \leq 2 \pi \}$ where
\[
k(\theta) = \left\{ 
\begin{array}{ll}
(\underline{k}(\theta),1) & \mbox{if } -\pi < \theta \leq \pi,
  \\ (\underline{k}(\theta),-1) & \mbox{if } -2\pi < \theta \leq -\pi
  \mbox{ or } \pi < \theta \leq 2 \pi.
\end{array}
\right.
\]
Then $K_\infty$ is a maximal compact subgroup of $G(\bbR)$ and
$\pr(K_\infty) = \underline{K}_\infty$.  If $r$ is an odd integer,
then $k(\theta) \mapsto e^{i\frac{r}{2}\theta}$ defines a genuine
character of $K_{\infty}$

Let $A_{r/2}(4)$ denote the set of functions $\varphi$ in
$L^2(\SL_2(\bbQ) \backslash G(\bbA))$ satisfying the following
properties:
\begin{enumerate}
\item $\varphi(g k_1) = \varphi(g)$ for all $k_1 \in K_1(4) \prod_{p
  \neq 2, \infty} K_p$,

\item $\varphi(g k_0) = \gamma(k_0) \varphi(g)$ for all $k_0 \in K_0$
  in $G(\bbQ_2)$ where $\gamma(-1) = - i^r$,

\item $\varphi(g r(\theta)) = e^{i \frac{r}{2} \theta} \varphi(g)$,

\item as a function on $G(\bbR)$, $\varphi$ is smooth and satisfies
  $\triangle \varphi = - \frac{r}{4} \left( \frac{r}{4} - 1 \right)
  \varphi$ where $\triangle$ is the Casimir operator and

\item $\varphi$ is cuspidal, ie. $\int_{N(\bbQ)\backslash N(\bbA)}
  \varphi(x(u)g)du = 0$ for all $g \in G(\bbA)$.
\end{enumerate}

A basis of $A_{r/2}(4)$ arises from cuspidal automorphic
representations $\pi=\otimes_{v}\pi_{v}$ of $G(\bbA)$ such that
$\pi_{\infty}$ is a holomorphic discrete series representation with
the lowest weight $r/2$, $\pi_{p}$ is unramified for all $p\neq 2$,
and $\pi_{2}$ contains a $K_{1}(4)$-fixed vectors.  In particular,
$\pi_{2}^{\gamma}\neq 0$ for some central character $\gamma$.  Note
that $\gamma$ is determined by $r$. Indeed, since the local components
of $\bfs_\bbQ(h(-1))$ for $v\neq \infty, 2$ are contained in $K_{p}$,
$\varphi(1) = \varphi(\bfs_\bbQ(h(-1))) = \gamma(-1) e^{i \pi r/2}
\varphi(1)$ and we get $\gamma(-1) = - i^r$.

\smallskip 
Let $\calH$ be the complex upper half plane. For $\underline{g}
= \begin{pmatrix} a & b \\ c & d \end{pmatrix} \in \SL_2(\bbR)$, $g
= (\underline{g}, \epsilon) \in G(\bbR)$ and $z \in \calH$, we define
\[
g z = \underline{g} z = \frac{az + b}{cz + d}.
\]
We define a holomorphic function on $\calH$ by
\[
J(g,z) = J((\underline{g},\epsilon),z) := \epsilon \, (cz + d)^{1/2}.
\]
Here we choose $w^{1/2}$ such that $- \frac{\pi}{2} < \arg(w^{1/2})
\leq \frac{\pi}{2}$. We call $J(g,z)$ a {\it factor of automorphy}. By
Lemma 3.3 in \cite{G2}, it satisfies $J(g g',z) = J(g, g'z) J(g',z)$
for any two $g$ and $g'$ in $G(\bbR)$.  Define a congruence subgroup
$\Gamma_{0}(4)$ by
\[
\Gamma_{0}(4) := G(\bbR) \cap (\bfs_\bbQ(\SL_{2}(\bbQ)) \cdot K_{0}(4)
\cdot \prod_{p\neq 2} K_{p}).
\] 
Similarly, define $\Gamma_{1}(4)\subseteq \Gamma_{0}(4)$ by replacing
$K_{0}(4)$ with $K_{1}(4)$.   Let $S_{r/2}(\Gamma_0(4))$ and
$S_{r/2}(\Gamma_1(4))$ be the spaces of classical modular forms of
weight $r/2$.  By page 183 in \cite{Kb}, $S_{r/2}(\Gamma_0(4)) =
S_{r/2}(\Gamma_1(4))$. We will denote this space by $S_{r/2}(4)$.

By Proposition 3.1 in \cite{G2}, there is a bijection $Q : A_{r/2}(4)
\rightarrow S_{r/2}(4)$ which we will recall below: Given $\varphi \in
A_{r/2}(4)$, then
\[
(Q \varphi)(z) = \varphi(g_\infty) J(g_\infty,i)^r
\]
where $z = g_\infty i \in \calH$.  Conversely, given $f \in
S_{r/2}(4)$. Let $g \in G(\bbA)$. By Lemma 3.2 in \cite{G2}, $g =
g_\bbQ g_\infty k$ for some $g_\bbQ \in \bfs_\bbQ(\SL_2(\bbQ))$,
$g_\infty \in G(\bbR)$ and $k \in K_1(4) \prod_{p \neq 2, \infty}
K_p$. Then
\[
(Q^{-1}f)(g) = f(g_\infty(i)) J(g_\infty,i)^{-r}.
\]

Using the bijection $Q$, we define another bijection between the
spaces of operators
\[
\bfq : \End_\bbC(A_{r/2}(4)) \rightarrow \End_\bbC(S_{r/2}(4))
\]
by $\bfq(L) = Q L Q^{-1}$.  Since the Hecke algebra $H(\gamma)$
defined in Section \ref{S3} acts on $A_{r/2}(4)$ it is of interest to
reinterpret this action in terms of classical modular forms.

\begin{prop} \label{P24} 
Let $U_{1}$ and $T_{1}$ be the operators in the local Hecke algebra
$H(\gamma)$ where $\gamma(-1)=-i^{r}$. Recall that 
$\zeta=\frac{1-i^r}{\sqrt{2}}$. 
For $f(z) \in S_{r/2}(4)$, we
have
\begin{enumerate}[(i)]
\item ${\displaystyle (\bfq(U_{1})f)(z) = \overzeta \,  (2z)^{-r/2} f
  \left( -\frac{1}{4z} \right)}$ and 

\item
${\displaystyle (\bfq(T_1)f)(z) = 2^{-r/2} \sum_{u = 0}^3
  f\left( \frac{z+u}{4} \right)}$.
\end{enumerate}
\end{prop}

\begin{proof}
  (i) Suppose $\varphi = Q^{-1}(f) \in A_{r/2}(4)$. For every place
  $v$, let $w_v = w(2^{-1})$ be the element in $G(\bbQ_v)$ defined in
  Section \ref{S1}.  By Proposition \ref{P3}(iv),
\[ 
(U_1 \varphi)(g_\infty)  =  \int_{K_0w_2 K_0} U_1(k)
  \varphi(g_\infty k) dk = U_1(w_2) \varphi(g_\infty w_2) = \overzeta
  \varphi(g_\infty w_2).
\]
Next, consider $\underline{w}(2^{-1})$ in $\SL_2(\bbQ)$.  
By (2.30) in \cite{G2}, $s_{\bbQ}(\underline{w}(2^{-1}))=
\prod w_v$. Since $\varphi$ is left $\SL_2(\bbQ)$-invariant, and 
right $K_p$-invariant for $p\neq 2$,  
\[
\overzeta\varphi(g_{\infty} w_2)=
\overzeta\varphi(s_{\bbQ}(\underline{w}(2^{-1}))^{-1} g_{\infty} w_2) =   
   \overzeta  \varphi
  (\left(\prod_{v \neq 2} w_v^{-1} \right) g_\infty) = 
  \overzeta  \varphi(w_\infty^{-1} g_\infty).
\] 
Applying $Q$ to the above equation gives (i). Part (ii) is proved
analogously.
\end{proof}

\section{Kohnen's plus space} \label{S8}

Hecke eigenforms in $S_{r/2}(4)$ correspond to cuspidal automorphic
representations $\pi$ such that $\pi_{\infty}$ is a discrete series
representation of lowest weight $\frac{r}{2}$, $\pi_{p}$ is unramified
for all $p\neq 2$, and $\pi_{2}$ has $K_{1}(4)$-fixed vectors. In
particular, $\pi_{2}^{\gamma}\neq 0$ for the central character
$\gamma(-1) = -i^{r}$. If $\pi_{2}$ is a principal series
representation then $\pi_{2}^{\gamma}$ is 2-dimensional and therefore
the corresponding Hecke eigenspace in $S_{r/2}(4)$ is also
2-dimensional.  Kohnen's plus space is introduced to resolve this
ambiguity. In terms of the space of automorphic functions
$A_{r/2}(4)$, it is clear what to do.  Decompose
\[
A_{r/2}(4)=A_{r/2}^{+}(4)\oplus A_{r/2}^{-}(4)
\]
where $A_{r/2}^{+}(4)$ is the eigenspace of the local Hecke operator
$T_{w}$ with the eigenvalue $2$, while $A_{r/2}^{-}(4)$ is the
eigenspace with the eigenvalue $-1$.  Since a presence of the
eigenvalue $2$ for $T_{w}$ acting on $\pi_{2}$ eliminates a
possibility that $\pi_{2}$ is a Steinberg representation, we see that
there is a one to one correspondence between Hecke eigenforms in
$A_{r/2}^{+}(4)$ and cuspidal automorphic representations $\pi$ (as
above) such that $\pi_{2}$ is a $\gamma$-unramified
representation. The classical Kohnen plus space is (essentially)
$Q(A_{r/2}^{+}(4))$ as it will be explained in a moment.
Niwa defines two operators ${\mathbf T}_{4}$ and ${\mathbf W}_4$ on
$S_{r/2}(4)$ \cite{Ni}:
\[
(\bfT_{4}f)(z) = \frac{1}{4} \sum_{u=0}^3 f\left(\frac{z+u}{4}\right)
\mbox{ \ and \ } (\bfW_{4}f)(z) = (-2iz)^{-r/2}
f\left(-\frac{1}{4z}\right).
\]
Note that $\bfW_4^2=1$. Let $\kappa = \frac{r-1}{2}$. Niwa shows that
the operator $\bfW= (-1)^{\frac{(r^2-1)}{8}} 2^{1-\kappa} \bfW_{4}
\bfT_{4}$ on $S_{r/2}(4)$ satisfies\footnote{In Kohnen's paper
  \cite{Ko}, the operator is $\bfT_4 \bfW_4$ acting on the right, ie
  $\bfT_4$ acts first and $\bfW_4$ follows.}  the quadratic relation
\[
(\bfW+1)(\bfW-2) = 0.
\]
Kohnen defines $S_{r/2}^+(4)$ and $S_{r/2}^-(4)$ to be the eigenspaces
of $\bfW$ on $S_{r/2}(4)$ of eigenvalues $2$ and $-1$ respectively
\cite{Ko}.  Proposition \ref{P24} says that
\[
\begin{cases} 
\bfq (U_1)= (-1)^{\frac{r^2-1}{8}} \bfW_4 \\
\bfq(T_1)= 2^{\frac{3}{2}-\kappa} \bfT_4 
\end{cases} 
\] 
 where the sign $(-1)^{\frac{r^2-1}{8}}$ is the
quotient of $\overzeta=\frac{1+i^r}{\sqrt{2}}$ and $i^{\frac{r}{2}}=
\left(\frac{1+i}{\sqrt{2}}\right)^r$.
Since $T_{w}=\sqrt{2}^{-1}T_{1}U_{1}$ it follows that 
$\bfq(T_{w})$ and $\bfW$ are conjugates of each other
by $\bfW_{4}$. 
  Thus the Kohnen's plus space is simply a conjugate of
our space:
\[
Q(A^{+}_{r/2}(4))=\bfW_{4}(S^{+}_{r/2}(4)).
\]
Since $\bfW_{4}$ commutes with the classical Hecke operators ${\mathbf
  T}_{p^{2}}$ for $p\neq 2$, $Q(A^{+}_{r/2}(4))$ and $S^{+}_{r/2}(4)$
are isomorphic as $\bbC[{\mathbf T}_{3^{2}}, {\mathbf T}_{5^{2}},
\ldots ]$-modules.
  
  There is another description of $S^{+}_{r/2}(4)$ in terms of Fourier
  coefficients.  It consists of the cusp
  forms whose $n$-th Fourier coefficients vanishes whenever
  $(-1)^{\kappa} n \equiv 2, 3 \pmod{4}$.  Kohnen also defines a Hecke
  operator ${\mathbf T}_4^+$ which preserves $S_{r/2}^+(4)$ in the
  following way: For $f(z) = \sum_n a_n q^n \in S_{r/2}^+(4)$, we set
  $({\mathbf T}_4^+ f)(z) = \sum_n b_n q^n$ where the sum is taken
  over integers $n > 0$ and $(-1)^{\kappa} n \equiv 0, 1 \pmod{4}$,
  and
\[
b_n = a_{4n} + \left( \frac{(-1)^{\kappa} n}{2} \right) 2^{\kappa-1}
a_n + 2^{r-2} a_{n/4}.
\]
Here $a_{n/4} = 0$ if $n$ is not a multiple of 4. The large
parenthesis denotes the Legendre symbol.

We can now formulate and prove our main global results.

\begin{thm} \label{T25}
There is a one to one correspondence between Hecke eigenforms $f$ in
$S_{r/2}^+(4)$ and irreducible cuspidal automorphic representations
$\pi = \otimes_v \pi_v$ in $L^2(\SL_2(\bbQ) \backslash G(\bbA))$ such
that
\begin{enumerate}[(i)]
\item $\pi_\infty$ is the discrete series representation of $G(\bbR)$
  with the lowest weight $\frac{r}{2}$. 

\item $\pi_p$ is unramified for all odd primes $p$.

\item $\pi_2$ is $\gamma$-unramified where $\gamma(-1) = - i^r$. 

\item If ${\mathbf T}_4^+ f = \lambda_2 f$, then a $\gamma$-spherical
  vector in $\pi_2$ is an eigenvector for $Z = \frac{T_{1}}{2} +
  (\frac{T_1}{2})^{-1}$ with the eigenvalue $2^{1-\frac{r}{2}}
  \lambda_2$.
\end{enumerate}
Note that $\lambda_2$ determines the eigenvalue of $Z$ on a
$\gamma$-spherical vector which in turn determines $\pi_2$ uniquely by
Corollary \ref{C20}. 
\end{thm}

\begin{proof} 
  The first three statements are clear, since $Q^{-1}(\bfW_{4}f)$ is a
  Hecke eigenform in $A_{r/2}^{+}(4)$ which is contained in a cuspidal
  automorphic representation $\pi$ with these properties. It remains
  to show (iv). We need the following lemma.

\begin{lemma} \label{L26}
Let $f$ be in $S^{+}_{r/2}(4)$. Then $\bfT_{4}^{+} f =
2^{\frac{r}{2}-1} \bfq(Z)f$.
\end{lemma}

\begin{proof}
Recall that $T_1$ is invertible by Proposition~\ref{P8}. Hence, it
suffices to show that
\[
2^{2 - r/2}
\bfq(T_1) {\mathbf T}_4^+ = \bfq(T_1^2 + 4).
\] 
 If $f(z) = \sum_{n=1}^\infty a_n q^n \in S_{r/2}(4)$ then, by
 Proposition~\ref{P24}, $(\bfq(T_1) f)(z) = 2^{2-r/2}
 \sum_{n=0}^\infty a_{4n} q^n$.  Thus, if $f(z) \in S_{r/2}^+(4)$,
 then one computes
\[
2^{2-r/2} (\bfq(T_1) \bfT_{4}^+ f)(z)=(\bfq(T_1^2 + 4)f)(z)=\sum_n
(2^{4-r} a_{16n} + 4a_n)q^n.
\]
 This proves the lemma. 
\end{proof}

Now we can finish the proof of Theorem \ref{T25}. If $\bfT_{4}^{+}f =
\lambda_{p}f$ then Lemma \ref{L26} implies that $Q^{-1}(f)$ is an
eigenform for $Z$ with the eigenvalue $2^{1-\frac{r}{2}} \lambda_{2}$.
Since $\bfW_{4} = (-1)^{\frac{r^2-1}{8}} \bfq(U_{1})$ and $Z$ commutes
with $U_{1}$, $Q^{-1}(\bfW_{4}f)$ is also an eigenform for $Z$ with
the same eigenvalue. This completes the proof of Theorem \ref{T25}.
\end{proof}

If $f$ is a Hecke eigenform in $S^{+}_{r/2}(4)$ then by Theorem 1(ii)
in \cite{Ko} the corresponding Shimura lift $f'= \rSh(f)$ is a Hecke
eigenform in $S_{r-1}(\SL_2(\bbZ))$. Recall that $G' = \PGL_2$. There
is a one to one correspondence between Hecke eigenforms $f'$ in
$S_{r-1}(\SL_2(\bbZ))$ and irreducible cuspidal automorphic
representations $\pi' = \otimes_v \pi_v'$ in $L^2(G'(\bbQ) \backslash
G'(\bbA))$ such that $\pi'_{\infty}$ is a discrete series
representation with the lowest weight $r-1$ and $\pi_{p}'$ is
unramified for all primes~$p$.  See Proposition~3.1 in \cite{G1}.  We
recall the local Shimura lift $\rSh(\pi_2)$ in Proposition~\ref{P11}
of a $\gamma$-unramified representation $\pi_2$ of $G(\bbQ_2)$. The
following corollary gives a precise representation-theoretic
description of the Shimura correspondence at the place $p=2$.

\begin{cor} \label{C27}
  Let $f$ be a Hecke eigenform in $S_{r/2}^+(4)$.  
\begin{enumerate}[(i)]
\item Let $\pi = \otimes_v \pi_v$ be the cuspidal automorphic
  representation corresponding to $f$ in Theorem \ref{T25}. 

\item Let $\pi' = \otimes_v \pi_v'$ be the cuspidal automorphic
  representations of $L^2(G'(\bbQ) \backslash G'(\bbA))$ corresponding
  to the Hecke eigenform $f' = \rSh(f)$ in $S_{r-1}(\SL_2(\bbZ))$.
\end{enumerate}
Then $\rSh(\pi_2) = \pi_2'$.
\end{cor}

\begin{proof} 
If $\bfT_4^{+} f = \lambda_{2}f$ then by Theorem 1(ii) in \cite{Ko},
$\bfT_{2} f' = \lambda_{2}f'$ where $\bfT_{2}$ is the classical Hecke
operator action on $S_{r-1}(\SL_2(\bbZ))$. By Proposition 5.2.1 in
\cite{G1} one checks that $\pi'_{2}$ is indeed isomorphic to
$\rSh(\pi_{2})$.
\end{proof}


Let $\pi$ be a cuspidal automorphic representation of $G(\bbA)$ as in
Theorem \ref{T25} and $\pi'$ be the corresponding cuspidal automorphic
representation of $G'(\bbA)$ as in Corollary \ref{C27}. By the
Ramanujan conjecture, proved by Deligne, $\pi_2'=\rSh(\pi_2)$ is a
tempered irreducible unramified representation so $\pi_2' = I'(s)$ for
some $s \in i \bbR$. This implies that $\pi_2 = I(\gamma,s)$ by
Theorem \ref{T22} and Corollary \ref{C23}. In particular
$\pi_2^\gamma$ is an irreducible $H(\gamma)$-module of dimension~2. It
corresponds under $Q$ to a two dimensional subspace of $S_{r/2}(4)$
spanned by a line in $S^+_{r/2}(4)$ and a line in $S^-_{r/2}(4)$.

On the other hand, if $\pi_2 = St(\epsilon)$ is a Steinberg
representation of $G(\bbQ_2)$ (see the definition before Corollary
\ref{C18}), then $\pi$ corresponds under $Q$ to an Hecke eigenform in
$S^-_{r/2}(4)$.  More precisely, we have the following theorem:

\begin{thm} \label{T28}
  There is a one to one correspondence between Hecke eigenforms $f$ in
  $S_{r/2}^-(4)$ such that $\bfW_4 f=-\epsilon(-1)^{\frac{r^2-1}{8}}
  f$, for some $\epsilon = \pm 1$, and irreducible cuspidal
  automorphic representations $\pi = \otimes_v \pi_v$ in
  $L^2(\SL_2(\bbQ) \backslash G(\bbA))$ such that
\begin{enumerate}[(i)]
\item $\pi_\infty$ is the discrete series representation of $G(\bbR)$
  with the lowest weight $\frac{r}{2}$.

\item $\pi_p$ is unramified for all odd primes $p$.

\item $\pi_2$ is the Steinberg representation $St(\epsilon)$.
\end{enumerate}
\end{thm}

\begin{proof}
  Recall, by Corollary \ref{Steinberg} that $T_w$ and $U_1$ act on
  one-dimensional space $St(\epsilon)^{\gamma}$ by $-1$ and
  $-\epsilon$. The theorem now follows form Proposition \ref{P24} and
  the definition of $S^-_{r/2}(4)$.
\end{proof}

\subsection*{Acknowledgment} 
The first author would like to thank the hospitality of the University
of Utah while part of this paper was written.  The second author is
supported by an NSF grant DMS-0551846.


\begin{thebibliography}{99}

\bibitem[BZ]{BZ} I. N. Bernstein and A. V. Zelevinsky, {\em
  Representations of the group $\GL(n,F),$ where $F$ is a local
  non-Archimedean field.} Russian Math. Surveys {\bf 31} (1976),
  no. 3, 1-68.

\bibitem[Bo]{Bo} A. Borel, {\em Admissible representations of a
    semi-simple group over a local field with vectors fixed under an
    Iwahori subgroup.} Invent. Math. {\bf 35}, (1976), 233-259. 


\bibitem[Ca]{Ca} W. Casselman, { \em Introduction to the theory 
of admissible representations of $p$-adic reductive groups.} 
Preprint. 

\bibitem[G1]{G1} S. Gelbart, {\em Automorphic forms on adele groups.}
  Annals of Mathematics Studies, No. 83. Princeton University Press,
  Princeton, N.J. (1975).

\bibitem[G2]{G2} S. Gelbart, {\em Weil's representation and the
  spectrum of the metaplectic group}. Lecture notes in Math.
  {\bf 530} Springer-Verlag (1976).


\bibitem[Kb]{Kb} N. Koblitz, {\em Introduction to elliptic curves and
  modular forms.} Graduate Texts in Mathematics, 97. Springer-Verlag,
  New York, (1984).

\bibitem[Ko]{Ko} W. Kohnen, {\em Modular forms of half-integral
  weight on $\Gamma_0(4)$.}  Math. Ann. {\bf 248} (1980), no. 3,
  249-266.

\bibitem[LS]{LS} H. Y. Loke and G. Savin, {\em Modular forms on
  nonlinear double covers of algebraic groups.} preprint.

\bibitem[Ma]{Ma} H. Matsumoto, {\em Analyse harmonique dans les
  syst\`{e}mes de Tits bornologiques de type affine.} Lecture
  Notes in Math. {\bf 590}. Springer-Verlag, Berlin-New York,
  1977.

\bibitem[Ni]{Ni} S. Niwa, {\em On Shimura's trace formula.}  Nagoya
  Math. J. {\bf 66} (1977), 183--202.

 
\bibitem[Sa]{Sa} G. Savin, {\em Lectures on representations of
    $p$-adic groups. Representations of real and $p$-adic groups},
  19--46, Lect.  Notes Ser. Inst. Math. Sci. Natl. Univ. Singap., 2,
  Singapore Univ.  Press, Singapore.


\bibitem[St]{St} M. Stein, {\em Surjective stability in dimension $0$
  for $K_{2}$ and related functors.}  Trans. Amer. Math. Soc.  {\bf
  178} (1973), 165--191.

\bibitem[Wa]{Wa} J.-L. Waldspurger, {\em Sur les coefficients de
  Fourier des formes modulaires de poids demi-entier.} 
  J. Math. Pures Appl. (9) {\bf 60} (1981), no. 4, 375-484.

\end{thebibliography}
\end{document}